\documentclass[]{article}

\usepackage[T2A]{fontenc}			
\usepackage[utf8]{inputenc}	
\usepackage[russian, english]{babel}	
\usepackage{amsmath,amsfonts,amssymb,amsthm,mathtools} 
\usepackage{indentfirst}
\usepackage{hyperref}
\usepackage{easyReview, soul, xcolor, todonotes}
\usepackage{comment}
\usepackage{gensymb}
\usepackage{tkz-euclide}
\usepackage{setspace,tikz,xcolor}
\usepackage{graphicx}
\usepackage[centertableaux]{ytableau}
\usetikzlibrary{arrows,matrix,math}
\tikzset{tab/.style={matrix of math nodes,column sep=-.35, row sep=-.35,text height=5pt,text width=5pt,align=center,inner sep=2,font=\footnotesize}}
\newif\iftikz
\tikztrue

\usepackage[style=ext-numeric, citestyle=numeric, backend=biber,articlein=false]{biblatex}
\newcommand{\GL}{GL}
\newtheorem{Th}{Theorem}
\newtheorem{cor}{Corollary}

\newcommand{\qbinom}[3]{\genfrac{[}{]}{0pt}{}{#1}{#2}_{#3}}

\newcommand{\qKrawtchouk}{$q$\nobreakdash-Krawtchouk }

\title{Skew Howe duality and q-Krawtchouk polynomial ensemble}
\author{Anton Nazarov$^{1,2}$, Pavel Nikitin$^{3}$, Daniil Sarafannikov$^{1,4}$\\
{\small$^{1}$Department of Physics, St. Petersburg State University,} \\
{\small  Ulyanovskaya 1, 198504 St.~Petersburg, Russia}\\
{\small$^{2}$email:antonnaz@gmail.com}\\
{\small$^{3}$ 
Laboratory of Representation Theory and Dynamical Systems,}\\
  {\small St. Petersburg Department
    of Steklov Mathematical Institute}\\
  {\small of Russian Academy of Sciences}\\
  {\small 191023, Fontanka 27, St. Petersburg, Russia}\\
  {\small email: pnikitin0103@yahoo.co.uk}\\
{\small$^{4}$ email: sardanis123@gmail.com}
}

\begin{document}
\binoppenalty=99999
\relpenalty=99999
\maketitle

\begin{abstract}
  We consider the decomposition into irreducible components of the
  exterior algebra
  $\bigwedge\left(\mathbb{C}^{n}\otimes
    \left(\mathbb{C}^{k}\right)^{*}\right)$ regarded as a
  $\GL_{n}\times\GL_{k}$ module. Irreducible $\GL_{n}\times\GL_{k}$
  representations are parameterized by pairs of Young diagrams
  $(\lambda,\bar{\lambda}')$, where $\bar{\lambda}'$ is the complement
  conjugate diagram to $\lambda$ inside the $n\times k$ rectangle. We
  set the probability of a diagram as a normalized specialization of
  the character for the corresponding irreducible component. For the
  principal specialization we get the probability that is equal to the
  ratio of the $q$-dimension for the irreducible component over the
  $q$-dimension of the exterior algebra. We demonstrate that this
  probability distribution can be described by the \qKrawtchouk
  polynomial ensemble. We derive the limit shape and prove the central
  limit theorem for the fluctuations in the limit when $n,k$ tend to
  infinity and $q$ tends to one at comparable rates.
\end{abstract}

\section*{Introduction and main results}
Various dualities play major role in asymptotic representation theory.
In particular, the Schur--Weyl duality between $\GL_{n}$ and
$\mathcal{S}_{k}$ was used by S. Kerov to study the distribution of
tensors by symmetry types \cite{kerov1986asymptotic}. If $n,k$ tend to
infinity with the same rate the limit shape of Young diagrams in the
decomposition of $\left(\mathbb{C}^{n}\right)^{\otimes k}$ into
irreducible $\GL_{n}$-modules coincides with the famous
Vershik--Kerov--Logan--Shepp limit shape \cite{vershik1977asymptotics,
  logan1977variational}. This is not the case if $k\sim n^{2}$, as
was demonstrated by P. Biane \cite{biane2001approximate}. The group
$\mathcal{S}_{k}$ is the Weyl group of $\GL_{k}$ therefore the Schur--Weyl duality
leads to the $(\GL_{n},\GL_{k})$ Howe duality
\cite{howe1989remarks,howe1995perspectives}, that is the decomposition
of the symmetric algebra
\begin{equation}
  \label{eq:1}
  S\left(\mathbb{C}^{n}\otimes \mathbb{C}^{k}\right)\cong
  \bigoplus_{\ell(\lambda)\leq\min(n,k)}V_{\GL_{n}}(\lambda)\otimes V_{\GL_{k}}(\lambda)
\end{equation}
into the multiplicity-free and sum of irreducible
$\GL_{n}\times \GL_{k}$ modules where the diagrams $\lambda$ have at
most $\min(n,k)$ rows. Restrict the decomposition (\ref{eq:1}) to the
diagrams of at most $m$ columns and consider the probability of a
diagram to be proportional to the dimension of an irreducible
component. Then this probability measure is the same as the measure on
the main diagonal of lozenge tilings of the hexagon with the sides
$(m,n,k,m,n,k)$ induced by the uniform measure \cite{cohn1998shape}.
The decomposition (\ref{eq:1}) is also related to celebrated Schur
measures \cite{okounkov2001infinite, Okounkov00}. Skew
$(\GL_{n},\GL_{k})$ Howe duality, that is the multiplicity-free
decomposition
\begin{equation}
  \label{eq:2}
  \bigwedge\left(\mathbb{C}^n \otimes (\mathbb{C}^k)^*\right) \cong
  \bigoplus_\lambda  V_{\GL_{n}}(\lambda)\otimes V_{\GL_{k}}(\bar{\lambda}'),
\end{equation}
is relatively less studied from the probabilistic point of view. The
measure on the diagrams $\lambda$ of size $m$ introduced as the ratio
of the dimension of the corresponding irreducible
$\GL_{n}\times\GL_{k}$ modules to the dimension of the $m$-th exterior
power was considered in \cite{panova2018skew}. Nevertheless the
relation between the measure that is given by the ratio of the
dimension of the irreducible module to the dimension of the whole
exterior algebra and the Krawtchouk polynomial ensemble does not
appear to be widely known before the paper \cite{nazarov2021skew}.

The decomposition (\ref{eq:2}) in terms of the characters is an
alternative form of the dual Cauchy identity for Schur polynomials
\cite{macdonald1998symmetric}. Therefore we can introduce the probability
measure as the ratio of characters 
\begin{equation}
  \label{eq:3}
  \mu_{n,k}(\lambda|\{x_{i}\}_{i=1}^{n},\{y_{j}\}_{j=1}^{k})=
  \frac{s_{\lambda}(x_{1},\dots,x_{n})s_{\bar{\lambda}'}(y_{1},\dots,y_{k})}
  {\prod_{i=1}^{n}\prod_{j=1}^{k}(x_{i}+y_{j})}.
\end{equation}
This measure (up to a change of $\bar{\lambda}'\to\lambda'$ and
$y_{j}\to 1/y_{j}$) was considered in \cite{GTW01}, but the limit shapes were not
explicitly discussed there. We consider the principal specialization
of characters of the form
$x_{i}=q^{i-1}, y_{j}=q^{j-1}$ and the specialization
$x_{i}=q^{i-1}, y_{j}=q^{1-j}$. As was demonstrated in \cite{nazarov2021skew},
the measures then take the form
\begin{equation}
  \label{eq:4}
  \mu_{n,k}(\lambda|q)=\frac{q^{||\lambda||}\dim_{q}\left(V_{\GL_{n}}(\lambda)\right)\cdot
    q^{||\bar{\lambda}'||}\dim_{q}\left(V_{\GL_{k}}(\bar{\lambda}')\right)}
  {\prod_{i=1}^{n}\prod_{j=1}^{k}(q^{i-1}+q^{j-1})}, 
\end{equation}
and
\begin{equation}
  \label{eq:5}
   \mu_{n,k}(\lambda|q,q^{-1})=\frac{q^{||\lambda||}\dim_{q}\left(V_{\GL_{n}}(\lambda)\right)\cdot
    q^{-||\bar{\lambda}'||}\dim_{1/q}\left(V_{\GL_{k}}(\bar{\lambda}')\right)}
  {\prod_{i=1}^{n}\prod_{j=1}^{k}(q^{i-1}+q^{1-j})},
\end{equation}
for each specialization respectively, where
$||\lambda|| = \sum_{i=1}^{n}(i-1)\lambda_i$ and $q$-dimension
$\mathrm{dim}_{q}$ will be defined in Section
\ref{sec:q-krawtch-polyn}. Our first result is the relation of these
measures to \qKrawtchouk ensembles. We recall that normalized \qKrawtchouk orthogonal polynomials $\tilde{K}^{q}_{j}(q^{-x};p,N,q)$ satisfy the following orthogonality relations~\cite[Section 14.15]{koekoek2010hypergeometric}
  \begin{equation}
    \label{eq:7}
    \sum_{i=0}^{N}\qbinom{N}{i}{q}p^{-i}q^{\binom{i}{2}-iN}\tilde{K}^{q}_{j}(q^{-i};p,N;q)
    \tilde{K}^{q}_{l}(q^{-i};p,N;q)=\delta_{jl}. 
  \end{equation}

\begin{Th}[\qKrawtchouk ensemble]
  \label{thm:q-Krawtchouk}
  The probability measure (\ref{eq:4}) defines a \qKrawtchouk
  polynomial ensemble,
  \begin{equation}
    \label{eq:6}
    \mu_{n,k}(\lambda|q)=\det\left(\sqrt{W(a_{i}) W(a_{j})}
      \sum_{l=0}^{n-1}\tilde{K}^{q}_{l}(a_{i})\;
      \tilde{K}^{q}_{l}(a_{j})\right)_{i,j=1}^{n},  
  \end{equation}
  where $a_{i}=\lambda_{i}+n-i$ and
  $\tilde{K}^{q}_{j}(x)=\tilde{K}^{q}_{j}(q^{-x};p,N;q)$ are the
  normalized $q$\nobreakdash-Krawtchouk polynomials with $N=n+k-1$, $p=q^{1-2n}$, and 
  $$
  W(a_i)=q^{\binom{a_i}{2}+a_i(n-k)}\qbinom{n+k-1}{a_i}{q}.
  $$
  The measure (\ref{eq:5}) defines a \qKrawtchouk polynomial ensemble for the
  polynomials ${\tilde{K}^{q}_{j}(q^{-a_{i}};q^{2-2n-k},n+k-1;q)}$.
\end{Th}
Next we describe the limit behavior of the correlation kernels of
\qKrawtchouk ensembles as $n,k\to\infty$, $q\to 1$ with compatible
rates. We prove the convergence of the determinantal point ensembles
to the limit determinantal random point process by the method of
Borodin and Olshaski~\cite{borodin2007asymptotics}, which uses the
spectral theory of self-adjoint operators on Hilbert space to
establish the pointwise convergence of the correlation kernels.

Our derivation of the limit correlation kernel for the \qKrawtchouk
ensemble is similar to the proof for the Charlie and Krawtchouk
ensembles in \cite{borodin2007asymptotics} and to the proof for the
Hahn ensemble in \cite{gorin2008nonintersecting}, therefore we present
only an outline of the proof.
\begin{Th}[Limit correlation kernel]
  \label{thm:limit-shape}
  As $ n,k \to \infty$ and $ q \to 1 $ in such a way that
  $ q=1-\frac{\gamma} {n}$ and $\frac{k}{n}\to c$, and the variables $a,b$
  are defined as $a=nt+u$, $b=nt+v$, and $t,u,v$ are finite, the correlation
  kernels
  $$\mathbf{K}_{n}(a,b)=\sqrt{W(a) W(b)}
  \sum_{l=0}^{n-1}\tilde{K}^{q}_{l}(a)\; \tilde{K}^{q}_{l}(b)$$
  converge to the discrete sine kernel
  \begin{equation}
    \label{eq:51}
    \lim_{n\to\infty}
    \mathbf{K}_{n}(nt+u,nt+v)=\mathbf{K}^{sine}_{\varphi}(u,v)=\frac{\sin(\varphi(u-v))}{\pi(u-v)}, 
  \end{equation}
  where for the measure \eqref{eq:4} we have
  \begin{equation}
    \label{eq:52}
    \varphi=\arccos\left(\mathrm{sgn}(-\gamma)\dfrac{e^{\gamma-\frac{\gamma  t}{2}}}{2}
      \dfrac{1-e^{\gamma(c-1)}}{\sqrt{(1-e^{\gamma t})(1-e^{\gamma(c+1-t)})}}\right),
  \end{equation}
 and for the measure \eqref{eq:5} we have
  \begin{equation}
    \label{eq:53}
    \varphi=\arccos\left(\mathrm{sgn}(-\gamma)
      \frac{e^{\frac{\gamma}{2}(t-c)}}{2}
      \dfrac{1-e^{\gamma c}-e^{\gamma (c-t)}+e^{\gamma (c+1-t)}}
      {\sqrt{(1-e^{\gamma t})(1-e^{\gamma(c+1-t)})}}\right).
  \end{equation}

\end{Th}

%%%%%
%%
%%  Namely, for all $ \varepsilon>0 $ we have
%%  \begin{equation}
%%    \mathbb{P}(d_{Q}(f_n,f)>\varepsilon)\rightarrow 0,\ \text{as}\ n\to \infty,
%%  \end{equation}
%%  where the distance $d_{Q}$ is defined by the quadratic functional
%%  \begin{equation}
%%    \label{eq:8}
%%    Q[f]=\int\int f'(x) f'(y)
%%    \left(\ln\left(\left(\frac{1-e^{-\gamma|x-y|}}{\gamma}\right)^{-1}\right)+
%%      \dfrac{\gamma}{2}(x+y-|x-y|)\right) \dx \dy 
%%  \end{equation}
%%  as $d_{Q}(f_{1},f_{2})=(Q[f_{1}-f_{2}])^{\frac{1}{2}}$.
%%
%%
%%    and similarly the distance is defined by a quadratic functional
%%  \begin{equation}
%%    \label{eq:37}
%%    \tilde{Q}[f]=??  %% Даниил, вы можете это вычислить?
%%  \end{equation}
%%%%%%

We also describe the global fluctuations around the average. The
change of diagram coordinates from $\{\lambda_{i}\}_{i=1}^{n}$ to
$\{a_{i}\}_{i=1}^{n}$ corresponds to the $45^{\degree}$ rotation and
then we scale the diagram by the factor $\frac{1}{n}$, switching to
coordinates $x_{i}=\frac{a_{i}}{n}$. Then the upper boundary of the
rotated and scaled diagram defines a continuous piecewise-linear
function $f_{n}\in C([0,c+1])$.
\begin{equation*}
    \begin{tikzpicture}[baseline=10, scale=0.3]
      \newcount\n
      \n=5;
      \draw[gray!40, very thin] (0,0) grid (12,8);  
      % \draw[gray!40, very thin] (\n,0) -- (\n+\k,\k) -- (\k,\k+\n) -- (0,\n);
      % \draw[gray!40, very thin] (0,0) grid (\n+\k,\n+\k);
      \foreach \y [count=\i from 0] in {7,4,3,3,1} {
        \draw[black, thin] (\n-\i,\i) -- (\n+\y-\i,\i+\y)--(\n-1+\y-\i,\i+1+\y)--
        (\n-1-\i,\i+1);
        \tikzmath{
          integer \k;
          \k=\n-\i;
          integer \ii;
          \ii=\i+1;
        }
        \draw[black, thick] (0,5)--(1,6)--(2,5)--(4,7)--(6,5)--(7,6)--(8,5)--(11,8)--(12,7);
        \draw[blue, thick, dashed] (\n-1+\y-\i,\i+1+\y) node[anchor= west, black] {$\lambda_{\ii}$} -- (\n-1+\y-\i,0) node[anchor=south west, black] {$a_{\ii}$} node{\textbullet};
        \foreach \len in {0,...,\y} {
          \draw[black,thin] (\n-\i+\len,\i+\len) -- (\n-1-\i+\len,\i+1+\len);
        }
      }
      \draw[black, thin] (\n,0)--(0,\n);
      \tkzInit[xmin=0,xmax=13,ymin=0,ymax=8.5]
      % \tkzGrid[sub,color=gray, very thin, subxstep=2,subystep=2]
      % \tkzAxeXY[thick]
      % \tkzLabelX[orig=true,label options={font=\tiny},step=2]
      % \tkzLabelY[orig=true,label options={font=\small},step=2]
      \tkzLabelX[orig=false,step=2]
      \tkzLabelY[orig=false,step=2]
      \tkzDrawX
      \tkzDrawY
    \end{tikzpicture}                              \begin{array}[b]{l}
                             \mbox{Diagram}\,\lambda=(7,4,3,3,1)\;\mbox{is rotated,}\\
                             \mbox{thick black line is upper boundary}\, f_{n},\\
                             \mbox{row lengths}\;\{\lambda_{i}\}\;
                             \mbox{correspond to the}\\
                             \mbox{point positions}\;\{a_{i}\}.
                           \end{array}
\end{equation*}

\begin{Th}[Central limit theorem]
  \label{central-limit-theorem}
  Consider random point processes $\left\{x_{i}=\frac{a_{i}}{n}\right\}_{i=1}^{n}$ corresponding to the probability distributions (\ref{eq:4}) or (\ref{eq:5}). Consider $\rho$ given by the formula \eqref{eq:limit-shape} or \eqref{eq:limit-shape-2} as a function of $e^{-\gamma t}$ and denote its support by $[b-2a,b+2a]$. Then for a linear statistics $X_{f}^{(n)}=\sum_{i=1}^{n} f(e^{-\gamma x_{i}})$, where $f\in C^{1}([b-2a,b+2a])$, we have
  \begin{equation}
    \label{eq:clt}
    X_{f}^{(n)}-\mathbb{E}X_{f}^{(n)}\to \mathcal{N}\left(0,\sum_{l\geq 1}l |\widehat{f}_{l}|^{2}\right),
  \end{equation}
  in distribution, as $n,k\to\infty$ with $c=\lim\frac{k}{n}$, where
  the Fourier coefficients $\widehat{f}_{l}$ are defined as
  \begin{equation}
    \label{eq:107}
    \widehat{f}_{l}=\frac{1}{2\pi}\int_{0}^{2\pi}
    f\left(2a\cos \theta+b\right)e^{-il\theta}\mathrm{d\theta}.
  \end{equation}
  The values $a,b$ are given by the formula (\ref{eq:47}) for the
  measure (\ref{eq:4}) and by formula~(\ref{eq:49}) for the measure
  (\ref{eq:5})
  
\end{Th}

The first correlation function gives us the limit density of points,
which is then used to write the explicit expression for the limit
shape $f(x)$. We do not present a proof of the uniform convergence of
random functions $f_{n}(x)$ to $f(x)$, since it requires a lot of
technical details and therefore will be presented in a separate
publication. But one can combine Theorems \ref{thm:limit-shape} and
\ref{central-limit-theorem} to obtain the weak convergence to the
limit shape.

\begin{cor}
  The limit shape of the upper boundary $ f_n $ of a rotated and
  scaled random Young diagram $\lambda$ with respect to the
  probability measure \eqref{eq:4} is given by the formula
  \begin{equation}
    \label{eq:20}
    f(x)=1+\int_{0}^{x}(1-2 \rho(t)\mathrm{d t}),
  \end{equation}
  where the limit density $ \rho(t) $ is given by the formula:
  \begin{multline}
    \label{eq:limit-shape}
    \rho(t)=\lim_{u\to v}
    \mathbf{K}^{sine}_{\varphi}(u,v)
    =\frac{\varphi}{\pi}=\\
    =\dfrac{1}{\pi}
    \arccos\left(\mathrm{sgn}(-\gamma)\dfrac{e^{\gamma-\frac{\gamma  t}{2}}}{2}
      \dfrac{1-e^{\gamma(c-1)}}{\sqrt{(1-e^{\gamma t})(1-e^{\gamma(c+1-t)})}}\right)
  \end{multline}
  For the probability measure \eqref{eq:5} the limit density is given
  by the formula
  \begin{equation}
    \label{eq:limit-shape-2}
    \rho(t)=\frac{1}{\pi}
    \arccos\left(\mathrm{sgn}(-\gamma)
      \frac{e^{\frac{\gamma}{2}(t-c)}}{2}
      \dfrac{1-e^{\gamma c}-e^{\gamma (c-t)}+e^{\gamma (c+1-t)}}
      {\sqrt{(1-e^{\gamma t})(1-e^{\gamma(c+1-t)})}}
    \right).
  \end{equation}
\end{cor}
\begin{figure}[h]
  \centering
  \includegraphics[width=0.49\linewidth]{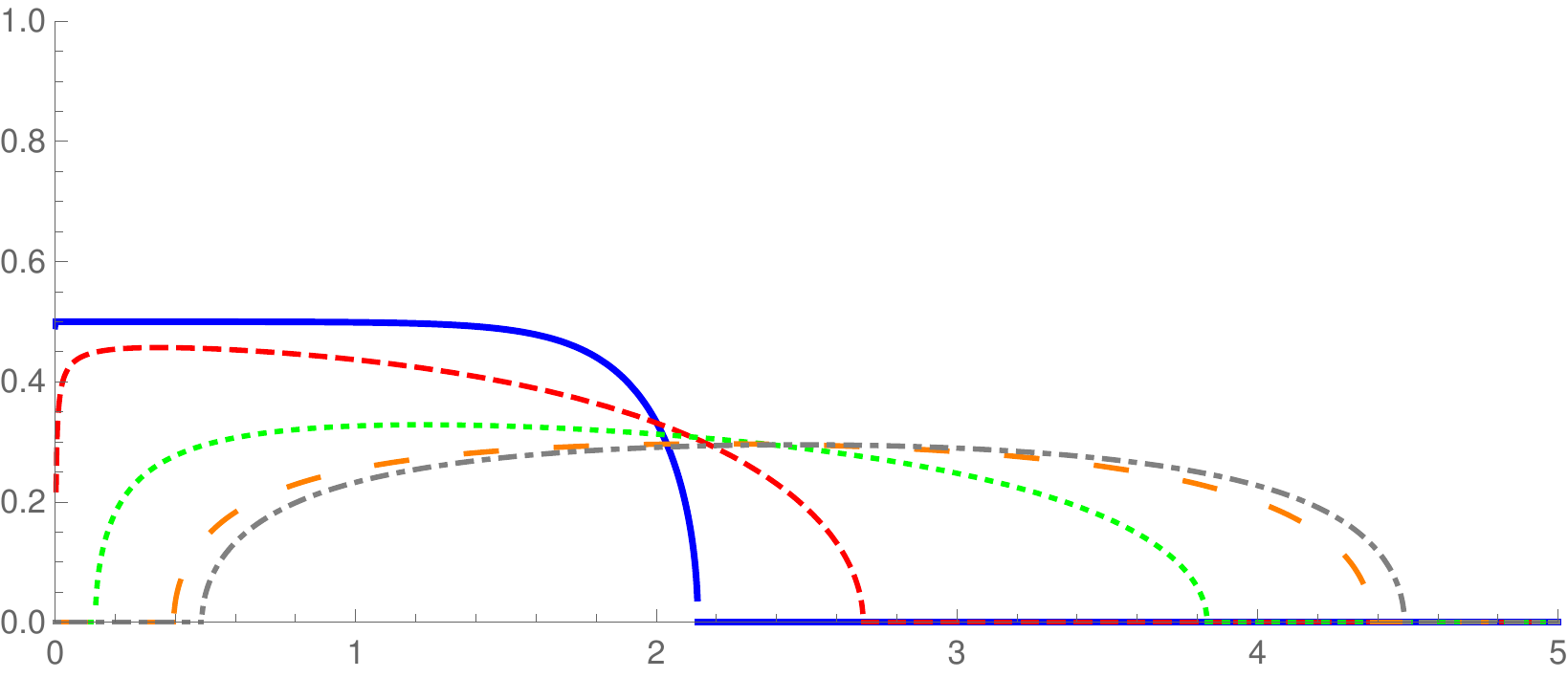}  
  \includegraphics[width=0.49\linewidth]{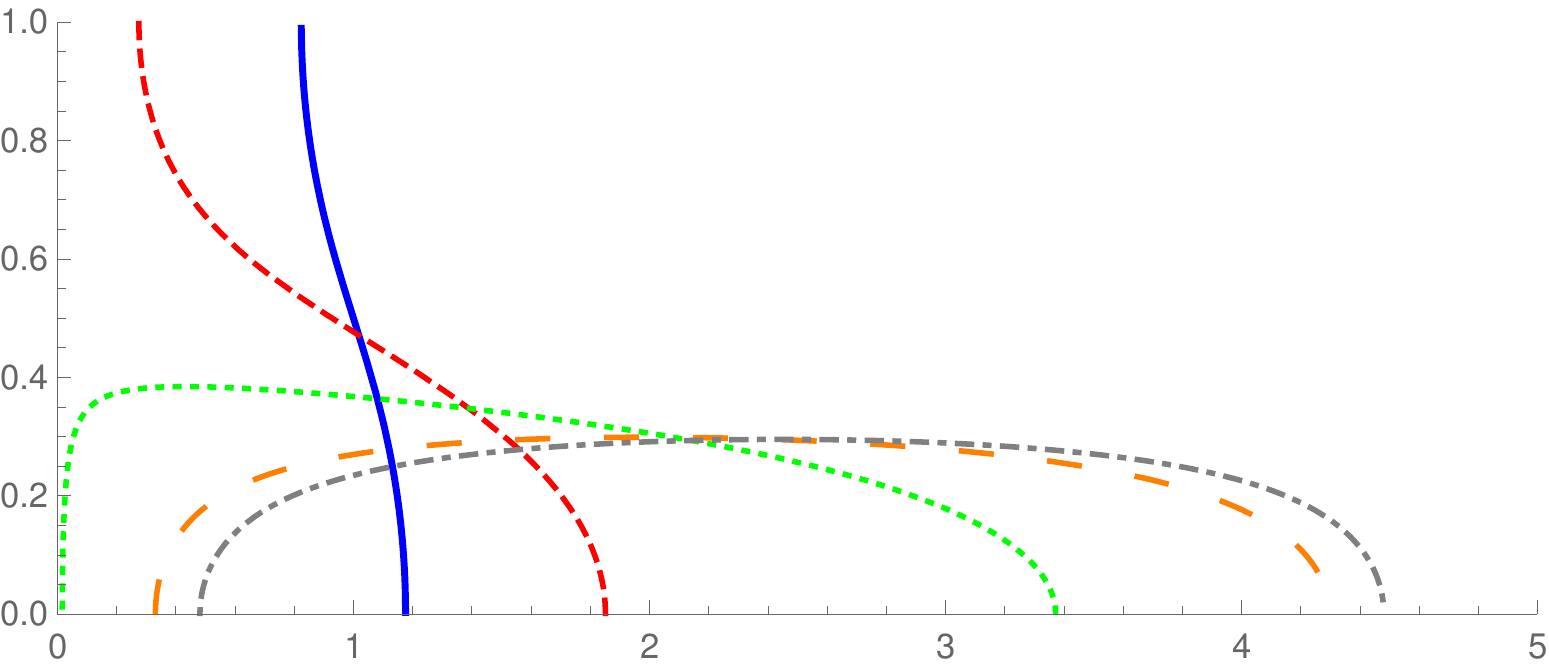}
  \caption{Plots of the limit densities \eqref{eq:limit-shape} (on the
    left) and \eqref{eq:limit-shape-2} (on the right) for $c=4$ and the
    values of $\gamma$: $-10$ (solid blue), $-2$ (dashed red), $-0.5$
    (dotted green), $-0.1$ (sparsely dashed orange) $-0.01$ (dot-dashed gray).  }
  \label{fig:limit-shape-densities}
\end{figure}

\begin{figure}[h]
  \centering
  \includegraphics[width=0.49\linewidth]{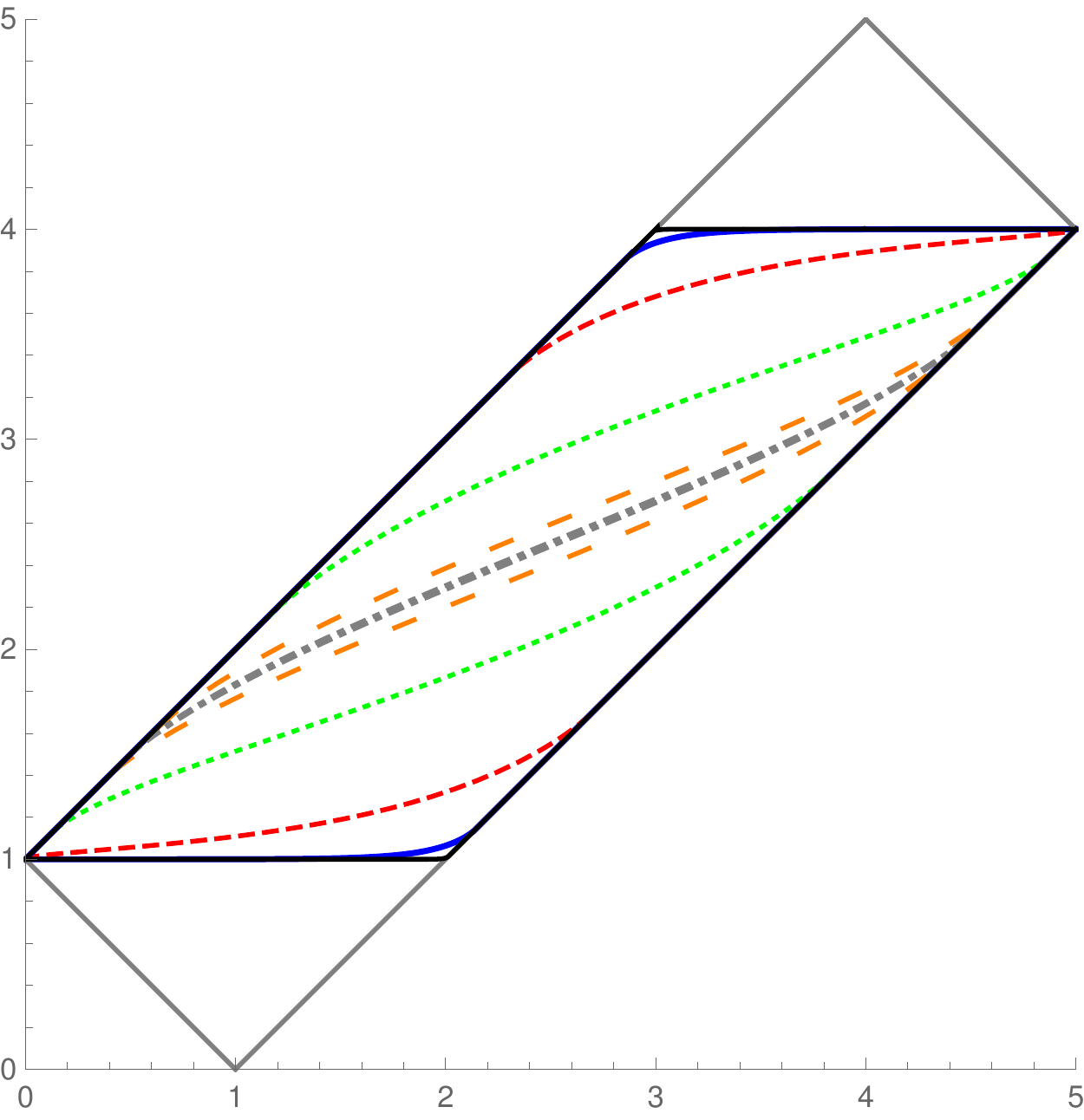}  
  \includegraphics[width=0.49\linewidth]{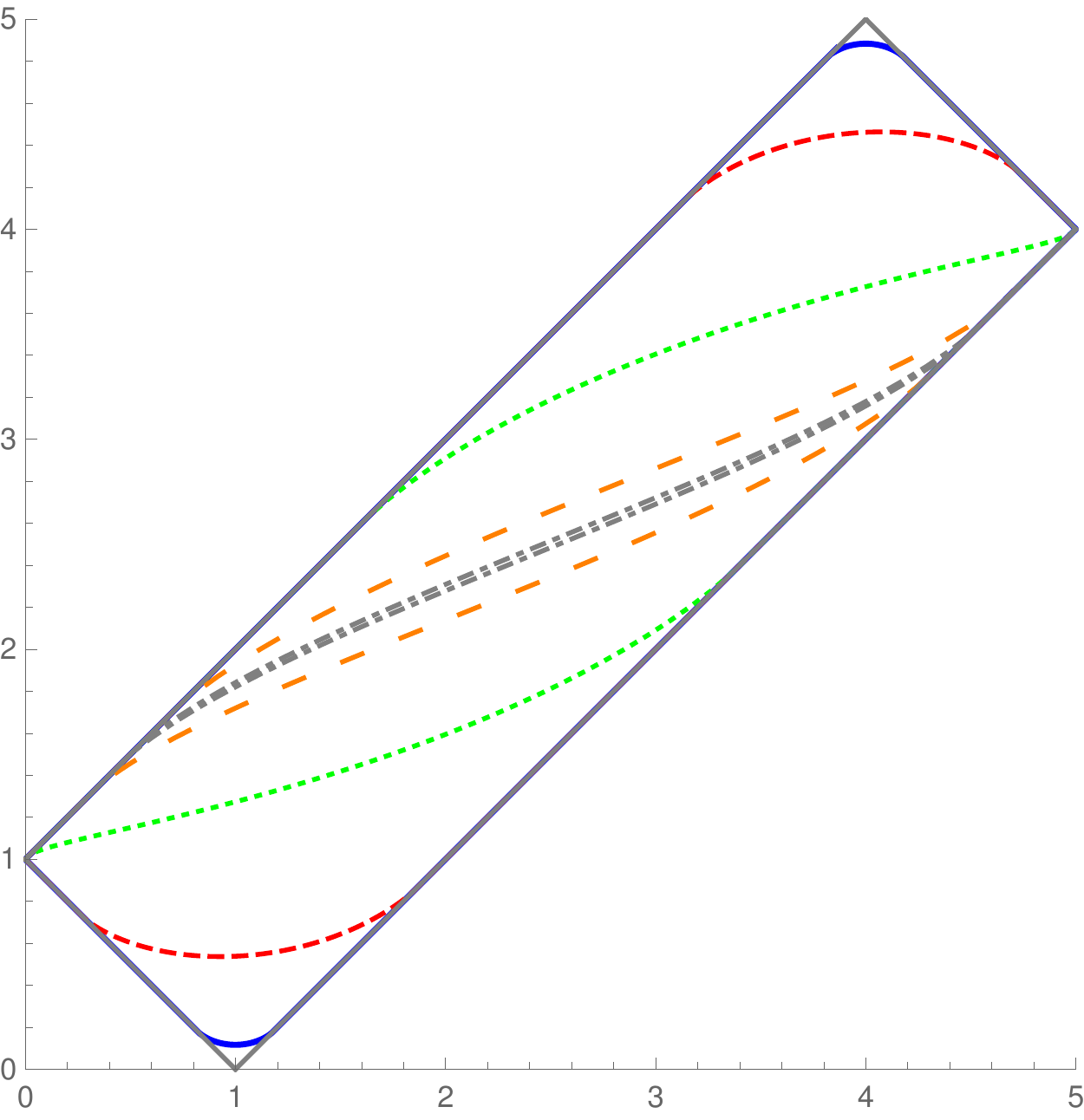}
  \caption{Plots of the limit shapes for Young diagrams corresponding
    to the densities \eqref{eq:limit-shape} (on the
    left) and \eqref{eq:limit-shape-2} (on the right) for $c=4$ and the
    values of $\gamma$ (bottom to top): $-10$ (solid blue), $-2$ (dashed red), $-0.5$
    (dotted green), $-0.1$ (sparsely dashed orange),  $-0.01$ (dot-dashed
    gray), $0.01$ (dot-dashed gray), $0.1$ (sparsely dashed orange), $0.5$
    (dotted green), $2$ (dashed red), $10$ (solid blue). Solid black
    lines on the left panel correspond to $\gamma=\pm\infty$
    ($q=\mathrm{const}$).  } %% Maybe add a phrase on q<1 q>1
  \label{fig:limit-shape-plots}
\end{figure}

Plots of the densities \eqref{eq:limit-shape} and
\eqref{eq:limit-shape-2} for $c=4$ are presented for various values of
$\gamma$ in Fig.~\ref{fig:limit-shape-densities}, and the corresponding
limit shapes are presented in Fig.~\ref{fig:limit-shape-plots}.
  
The paper is organized as follows. In Section
\ref{sec:q-krawtch-polyn} we use the explicit formulas for the
$q$-dimensions to prove Theorem \ref{thm:q-Krawtchouk}. In 
Section \ref{sec:limit-shape-young} we derive the limit shapes and
outline the proof of Theorem \ref{thm:limit-shape}. Then in Section
\ref{sec:fluctuations} we discuss the fluctuations and prove the central
limit theorem. We discuss some open questions in the conclusion.

\section{\qKrawtchouk polynomial ensemble}
\label{sec:q-krawtch-polyn}

In this section we prove Theorem \ref{thm:q-Krawtchouk}. We first recall the derivation of the explicit formula for the measures \eqref{eq:4}, \eqref{eq:5} from \cite[Theorem 4.6]{nazarov2021skew}  and then use it to demonstrate \eqref{eq:6}.

It is well known that the $\GL_{n}$ character, which is the Schur polynomial $s_{\lambda}(x_{1},\dots,x_{n})$, is given by the following sum over the semi-standard Young tableaux of the shape $\lambda$:
\begin{equation}
  \label{eq:9}
  s_{\lambda}(x_{1},\dots,x_{n})=\sum_{T\in SSYT(\lambda,n)} \prod_{i=1}^{n}
  x_{i}^{\#i's\ \text{in}\ T}.
\end{equation}
Define $q$-dimension of the irreducible $\GL_{n}$ representation as the
principal gradation (see \cite[\S10.10]{kac1990infinite}) that is the weighted
sum of the dimensions of weight subspaces:
\begin{equation}
  \label{eq:11}
  \dim_{q}\left(V_{\GL_{n}}(\lambda)\right)=\sum_{(u_{1},\dots,u_{n-1})\in
    \mathbb{Z}^{n-1}_{\geq 0}}q^{\sum_{i=1}^{n-1}u_{i}}\dim
  V(\lambda)_{\lambda-\sum_{i=1}^{n-1}u_{i}\alpha_{i}}, 
\end{equation}
where $\alpha_{1},\dots,\alpha_{n-1}$ are the simple roots of $\GL_{n}$ and we
identify the diagram $\lambda$ with the $\GL_{n}$ weight $\lambda$. We use the
notation
\begin{equation}
  \label{eq:q-number}
  [m]_{q}=\frac{1-q^{m}}{1-q}
\end{equation}
for the $q$-numbers, define $q$-factorials as products of $q$-numbers
and $q$-binomial coefficients as the ratio of $q$-factorials. The
formulas for \qKrawtchouk polynomials in
~\cite[Section 14.15]{koekoek2010hypergeometric} use
$q$-Pochhammer symbols, defined as
\begin{equation}
  \label{eq:45}
  (a;q)_{0}=1,\quad (a;q)_{m}=\prod_{i=1}^{m}(1-aq^{i-1}),\;
  m\in\mathbb{Z}_{+}.
\end{equation}

\begin{figure}[htb]
  \[
{\small\ytableaushort{11144,2246,33,44,6} \quad }
\begin{tikzpicture}[baseline=3em, scale=0.3]
  \draw[->, thick] (0,0) -- (11.5,0);
  \draw[->, thick] (0,0) -- (0,7.5);  
\draw[gray!40, very thin] (0,0) grid (11,7);
% Draw the paths
\draw[very thick, blue, line join=round] (5,1) -- (8,1) -- (8,4) -- (10,4) -- (10,7);
\draw[very thick, blue, line join=round] (4,1) -- (4,2) -- (6,2) -- (6,4) -- (7,4) -- (7,6) -- (8, 6) -- (8,7);
\draw[very thick, blue, line join=round] (3,1) -- (3,3) -- (5,3) -- (5,7);
\draw[very thick, blue, line join=round] (2,1) -- (2,4) -- (4,4) -- (4,7);
\draw[very thick, blue, line join=round] (1,1) -- (1,6) -- (2,6) -- (2,7);
\draw[very thick, blue, line join=round] (0,1) -- (0,7);
\draw(0-0.35,0) node[anchor=north] {$0$};
\foreach \i in {0,1,2,3,4,5} {
  \draw[fill=red, color=red] (\i, 1) circle (0.12);
%%  \draw (\i, 1) node[anchor=east] {$s_{\i}$};
}
\foreach \i in {0,2,4,5,8,10} {
  \draw[fill=red, color=red] (\i, 7) circle (0.12);
%%  \draw (\i, 1) node[anchor=east] {$s_{\i}$};
}
\end{tikzpicture}\quad
\begin{tikzpicture}[baseline=3em, scale=0.3]
  \draw[->, thick] (0,0) -- (14.5,0);
  \draw[->, thick] (0,0) -- (0,7.5);  
\draw[gray!40, very thin] (-2,0) grid (14,7);
% Draw the paths

\draw[very thick, blue, line join=round] (2,1) -- (1,2) -- (0,3) -- (-1,4) -- (-2,5) -- (-1,6) -- (-2,7);

\draw[very thick, blue, line join=round] (4,1) -- (3,2) -- (2,3) -- (1,4) -- (2,5) -- (3,6) -- (4,7);

\draw[very thick, blue, line join=round] (6,1) -- (5,2) -- (7,4) -- (6,5) -- (8,7);

\draw[very thick, blue, line join=round] (8,1) -- (11,4) -- (10,5) -- (11,6) -- (10,7);

\draw[very thick, blue, line join=round] (10,1) -- (13,4) -- (12,5) -- (14,7);

\draw(0-0.35,0) node[anchor=north] {$0$};
\foreach \i in {1,2,3,4,5} {
  \draw[fill=red, color=red] (2*\i, 1) circle (0.12);
%%  \draw (\i, 1) node[anchor=east] {$s_{\i}$};
}
\foreach \i in {-2,4,8,10,14} {
  \draw[fill=red, color=red] (\i, 7) circle (0.12);
%%  \draw (\i, 1) node[anchor=east] {$s_{\i}$};
}
\end{tikzpicture}
\]

  \caption{A $\GL_{6}$-Young tableau, its row reading ({\it left}) and its column readings ({\it right}).}
  \label{fig:row-col-reading}
\end{figure}
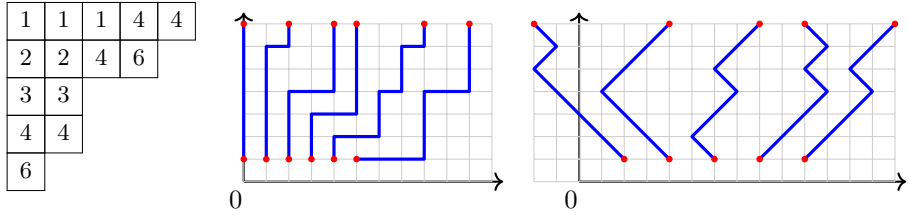
By the standard row-reading rule where the number of $i$-boxes in the
row $j$ corresponds to the number of horizontal steps on the level $i$
in the path number $j$ (see Fig.~\ref{fig:row-col-reading} ({\it
  left})), semistandard Young tableaux $SSYT(\lambda,n)$ are in
one-to-one correspondence with the configurations of $n$
non-intersecting paths that start at points
$(0,1),(1,1),\dots,(n-1,1)$ and end at points
$(a_{n},n),\dots,(a_{1},n)$. Horizontal steps on the level $i$ are
weighted by $x_{i}=q^{i-1}$ and vertical steps have weight 1. Using
the Lindstr\"om--Gessel--Viennot lemma~\cite{GV85,Lindstrom73} and the
recursion on the determinants it is easy to derive a well-known
$q$-analog of the Weyl dimension formula~\cite{kirillov2008introduction,BKW16}:
\begin{equation}
  \label{eq:10}
  s_{\lambda}(1,q,q^{2},\dots,q^{n-1})=q^{||\lambda||}\dim_{q}\left(V_{\GL_{n}}(\lambda)\right)=
  q^{||\lambda||}\prod_{i=1}^{n}\prod_{j=i+1}^{n}\frac{[a_{i}-a_{j}]_{q}}{[j-i]_{q}}.
\end{equation}
Similarly we get the formula for the $q$-dimension of the
$\GL_{k}$-representation $V_{\GL_{k}}(\bar{\lambda}')$, but we would
like to have it in terms of row lengths of $\lambda$, not
$\bar{\lambda}'$. Therefore we use the column reading for the bijection
between non-intersecting lattice paths and semistandard Young
tableaux. For the semistandard Young tableau $T\in SSYT(\mu,k)$ of at
most $n$ columns of lengths $\mu'_{1},\dots,\mu'_{n}$ the paths start
at points $(2,1),(4,1),\dots (2n,1)$ and go to the points
$(2i-2\mu'_{i}+k,k+1)$ for $i=1,\dots,n$. The paths consists of $k$
steps $(-1,1)$ or $(1,1)$ and the step number $j$ is $(-1,1)$ if $j$
is present in the column (see Fig.~\ref{fig:row-col-reading} ({\it
  right})). We weight the $(-1,1)$ steps by $q^{m/2}$ where $m$ is the
number of the left-leaning diagonal starting from $0$. We again can
use Lindstr\"om--Gessel--Viennot lemma and Dodgson condensation (see,
e.g.~\cite{BKW16}) for the determinants to obtain the formula
\begin{multline}
  \label{eq:12}
  s_{\bar{\lambda}'}(1,q,q^{2},\dots,q^{k-1})=
  q^{||\bar{\lambda}'||}\dim_{q}\left(V_{\GL_{k}}(\bar{\lambda}')\right)=\\
  = q^{||\bar{\lambda}'||}\prod_{1\leq i<j\leq n}^{n}[a_{i}-a_{j}]_{q}\cdot
  \prod_{l=1}^{n}\frac{[n+k-l]_{q}!}{[a_{l}]_{q}![n+k-1-a_{l}]_{q}!}.
\end{multline}

Now we can rewrite the measure \eqref{eq:4} as
\begin{equation}
  \label{eq:13}
  \mu_{n,k}(\lambda|q)=\frac{q^{\lVert\lambda\rVert+\lVert\bar\lambda^{'}\rVert}}
  {\prod_{i=1}^{n}\prod_{j=1}^{k}(q^{i-1}+q^{j-1})}
  \times\frac{\prod_{m=0}^{n-1}[k+m]_q!\prod_{1\leq i<j \leq n}[a_i-a_j]_q^2}
  {\prod_{i<j}[j-i]_q\prod_{i=1}^{n}[a_i]_q![n+k-1-a_i]_q!},
\end{equation}
where
\begin{equation}
  \label{eq:14}
  \lVert\bar\lambda^{'}\rVert =
  \sum_{i=1}^{n}\frac{(k-\lambda_i)(k-\lambda_i-1)}{2}= 
  \sum_{i=1}^{n}\binom{n+k-a_i-i}{2}.
\end{equation}
We rewrite the power of $q$ as
\begin{equation}
  \label{eq:15}
  q^{\lVert\lambda\rVert+\lVert\bar\lambda^{'}\rVert}=q^{\sum_{i=1}^{n}(i-1)(a_i-n+i)+\binom{n+k-a_i-i}{2}} 
  \propto q^{\sum_{i=1}^{n}\binom{a_i}{2}+2a_i(i-n)+a_i(n-k)},
\end{equation}
and $q$-analog of the Vandermonde determinant as
\begin{multline}
  \label{eq:16}
  \prod_{1\leq i<j \leq n}[a_i-a_j]_q^2 \propto
  \prod_{1\leq i<j \leq n}(1-q^{a_i-a_j})^2=\\
=  \prod_{1\leq i<j \leq n}q^{2a_i}(q^{-a_i}-q^{-a_j})^2
  =
  q^{\sum_{i=1}^{n}2a_i(n-i)}\prod_{i<j}(q^{-a_i}-q^{-a_j})^2,
\end{multline}
to write the measure in the form of a determinantal ensemble 
\begin{equation}
  \label{eq:17}
  \mu_{n,k}(\{a_i\},q)=C_{n,k,q}\prod_{i<j}(q^{-a_i}-q^{-a_j})^2\prod_{i=1}^{n}W(a_i),
\end{equation}
where
\begin{equation}
  \label{eq:18}
  W(a_i)=q^{\binom{a_i}{2}+a_i(n-k)}\qbinom{n+k-1}{a_i}{q}
\end{equation}
and
\begin{equation}
  \label{eq:19}
  C_{n,k,q} =
  \dfrac{q^{\frac{k n}{2}(n+k-2)}}
  {\prod_{i=1}^{n}\prod_{j=1}^{k}(q^{i-1}+q^{j-1})}
  \prod_{i=1}^{n}\dfrac{\left[k+i-1\right]_q!}
  {\left[i-1\right]_q! \left[n+k-1\right]_q!}\frac{1}{(1-q)^{\frac{n(n-1)}{2}}}.
\end{equation}
The weight $W(a_{i})$ coincides \cite[(14.15.2)][]{koekoek2010hypergeometric} with the weight for \qKrawtchouk polynomials $K_{m}^{q}(q^{-x};p,N;q)$ with the parameters $p=q^{1-2n}$ and $ N=n+k-1$. Therefore the equality \eqref{eq:6} is proven.

Similarly, the measure $\mu(\lambda|,q,q^{-1})$ is written explicitly as
\begin{equation}
  \label{eq:21}
  \mu_{n,k}(\lambda|q,q^{-1})=\hat{C}_{n,k,q}\prod_{i<j}(q^{-a_i}-q^{-a_j})^2\prod_{i=1}^{n}\hat{W}(a_i), 
\end{equation}
where
\begin{equation}
  \label{eq:22}
  \hat{W}(a_i) =q^{\binom{a_i}{2}+a_i(n-1)}\qbinom{n+k-1}{a_i}{q}
\end{equation}
and
\begin{equation}
  \label{eq:23}
  \hat{C}_{n,k,q}= 
  \dfrac{q^{\frac{n}{2}(n-1)(n+2k-2)}}
  {\prod_{i=1}^{n}\prod_{j=1}^{k}(q^{i-1}+q^{1-j})}
  \prod_{i=1}^{n}\dfrac{\left[k+i-1\right]_{q}!}
  {\left[i-1\right]_q! \left[n+k-1\right]_{q}!}\frac{1}{(1-q)^{\frac{n(n-1)}{2}}}.
\end{equation}
Taking $p=q^{2-2n-k}$ and $N=n+k-1$ we conclude the proof of Theorem
\ref{thm:q-Krawtchouk}.

\section{Correlation kernels and limit density}
\label{sec:limit-shape-young}
In this section we outline the proof of Theorem \ref{thm:limit-shape}. As was shown in the previous section, the upper boundary of a rotated random Young diagram corresponds to a point configuration. Therefore to derive the limit shape it is sufficient to find the limit density of the points. We use the $q$-difference equation for \qKrawtchouk polynomials to derive the limit density by the method of Borodin and Olshanski~\cite{borodin2007asymptotics}. The limit density is given by the
discrete sine-kernel as one expects from its universal properties~\cite{baik2007discrete}.

%%Then we consider the variational
%%problem, demonstrate non-negativity of the quadratic functional and
%%prove the uniform convergence of the random Young diagrams to the
%%limit shape with respect to the norm defined by the quadratic
%%functional.
%%
For any $n$-point discrete determinantal polynomial ensemble $\mathcal{P}^{(n)}$
with the weight function $W^{(n)}(x)$ and normalized orthogonal
polynomials $p^{(n)}_{i}(x)$ defined on a finite lattice
$\{x^{(n)}_{0},\dots,x^{(n)}_{L}\}$ , $m$-point correlation function
can be written as a determinant
\begin{equation}
  \label{eq:24}
  \rho^{(n)}_{m}(x_1,\ldots,x_m) = \det[\mathbf{K}_n(x_i,x_j)]_{1 \le i,j \le m},
\end{equation}
where $\mathbf{K}_n(x_i,x_j) $ is the correlation kernel defined by the formula
\begin{equation}
  \label{eq:25}
  \mathbf{K}_n(x,y) = \sum_{i=0}^{n-1} \sqrt{W^{(n)}(x)W^{(n)}(y)}\ p^{(n)}_i(x) p^{(n)}_i(y).
\end{equation}
The variables $x,y$ take values on the lattice
$\{x^{(n)}_{0},\dots,x^{(n)}_{L}\}$, therefore we can consider the
polynomials and the weights $W^{(n)}$ to be the functions of the
corresponding integer index $p^{(n)}_{i}(x_{l})=p^{(n)}_{i}(l)$. Then
the functions $\{\sqrt{W^{(n)}(x)}\ p^{(n)}_i(x)\}_{i=0}^L$ form an
orthonormal basis in the space $\ell^{2}$ on the finite set
$\{0,\dots,L\}$. The correlation kernel acts on this space by
projecting to the subspace spanned by the first $n$ states
$ \{\sqrt{W^{(n)}(x)}\ p^{(n)}_i(x)\}^{n-1}_{i=0} $. To prove the
convergence of the determinantal ensembles $\{\mathcal{P}^{(n)}\}$ as
$n\to\infty$ to a determinantal point process, it is sufficient to
demonstrate the pointwise convergence of the correlation kernels with
an appropriate scaling of the arguments
$\mathbf{K}_{n}(nt+x,nt+y)\xrightarrow[n\to\infty]{}
\mathbf{K}_{t}(x,y)$ \cite{borodin2007asymptotics,
  anderson2010introduction}. Assume that as $n\to\infty$, the lattice
also grows, so $L\to\infty$ as well and its points fill some interval.
Then the limit density of points in the point ensemble can be recovered
from the 1-point function $\rho(t\vert y)=\lim\limits_{x\to y}\mathbf{K}_{t}(x,y)$.

Consider the continuation of the functions $W^{(n)}, p^{(n)}_{i}$ to $\mathbb{Z}_{+}$ by assuming zero values at $x_{i}$ for ${i>L}$. Assume that for any $n$ there exists a bounded self-adjoint operator $D^{(n)}$ in
$\ell^{2}(\mathbb{Z}_{+})$ such that
$\{\sqrt{W^{(n)}(x)}\ p^{(n)}_i(x)\}_{i=0}^L$ are its eigenfunctions and
assume that there exits a limit $D^{(n)}\xrightarrow[n\to\infty]{} D$ in strong resolvent sense, where $D$ is a bounded self-adjoint operator on $\ell^{2}(\mathbb{Z}_{+})$ with a simple continuous spectrum $[\alpha,\beta]$. Then the theorem VIII.24 in \cite{reed1981} implies the convergence of the corresponding spectral projections.

The limit correlation kernel is given by a spectral projection to a subinterval of $[\alpha,\beta]$. Moreover, Hilbert space $\ell^{2}(\mathbb{Z}_{+})$ is isomorphic to $L^{2}([\alpha,\beta],d\nu)$, where $d\nu$ is the spectral measure on $[\alpha,\beta]$. The operator $D$ on
$L^{2}([\alpha,\beta],d\nu)$ becomes a multiplication operator and
spectral projection is given by the characteristic function. Taking
the Fourier transform from $L^{2}$ to $\ell^{2}$, we can recover the
limit correlation kernel $\mathbf{K}(x,y)$. 

In our case we consider the ensemble of \qKrawtchouk polynomials,
defined on the lattice $q^{-a}=e^{\gamma\frac{a}{n}}$,
$a=0,\dots,n+k-1$ with the weight \eqref{eq:18} or \eqref{eq:22}.
Since number of lattice points grows with $n$, we need to rescale our
problem and consider random variables $x=\frac{a}{n}$ that take
values on the interval $[0,c+1)$. Recall that $q$-difference equation for the \qKrawtchouk polynomials $K_{m}^{q}(q^{-a};p,N;q)=K_{m}^{q}(a)$ is given by \cite[(14.15.5)][]{koekoek2010hypergeometric}
\begin{equation}
  \label{eq:26}
  A(m)K_m^q(a)=B(a)K_m^q(a+1)-(B(a)+C(a))K_m^q(a)+C(a)K_m^q(a-1),
\end{equation}
where we have omitted some arguments of the coefficients $A,B,C$ for brevity:
\begin{align}
  \label{eq:27}
  A(m)&=  A(m,q)   = q^{-m}(1-q^m)(1+pq^m),\\
  B(a)&=  B(q,a,N) = 1-q^{a-N},\\
  C(a)&=  C(q,a)   = -p(1-q^a).
\end{align}
Rewriting for the functions
$ \kappa_m(a)=\sqrt{W(a)}\tilde{K}_{m}^{q}(a)$ and canceling the
normalization constant, we obtain
\begin{equation}
  \label{eq:28}
  \dfrac{A(m)}{\sqrt{W(a)}}\kappa_m(a)=
  \dfrac{B(a)}{\sqrt{W(a+1)}}\kappa_m(a+1)-\dfrac{(B(a)+C(a))}{\sqrt{W(a)}}
  \kappa_m(a)+\dfrac{C(a)}{\sqrt{W(a-1)}}\kappa_m(a-1).
\end{equation} 
Then we move some terms to the other side and multiply both sides by $ \sqrt{W(a)} $ to get
\begin{multline}
  \label{eq:29}
  B(a)\sqrt{\dfrac{W(a)}{W(a+1)}}\kappa_m(a+1)+
  C(a)\sqrt{\dfrac{W(a)}{W(a-1)}}\kappa_m(a-1)=\\
  =(A(m)+B(a)+C(a))\kappa_m(a). 
\end{multline}
If we express $W(a+1)$ as a product of $W(a)$ and the remaining term,
\begin{eqnarray}
  \label{eq:30}
  W(a+1)=W(a)q^{a+1-N}p^{-1}\dfrac{\left[N-a\right]_{q}}{\left[a+1\right]_{q}}=
  W(a)q^{a+1-N}p^{-1}\dfrac{1-q^{N-a}}{1-q^{a+1}},\\
  W(a-1)=W(a)q^{N-a}p\dfrac{1-q^a}{1-q^{N+1-a}},
\end{eqnarray}
it is easy to check that
\begin{equation}
  \label{eq:55}
  B(a)\frac{\sqrt{W(a)}}{\sqrt{W(a+1)}}=C(a+1)\frac{\sqrt{W(a+1)}}{\sqrt{W(a)}}=
  -\sqrt{pq^{a-N+1}(1-q^{a+1})(1-q^{N-a})}.   
\end{equation}

Now the left hand side of \eqref{eq:29} can be seen as an action of an
operator $D^{(n)}$ in $\ell^{2}(\mathbb{Z}_{+})$ on its eigenfunction
$\kappa_{m}(a)$. This action can be seen as a convolution with the
matrix $D^{(n)}(a,b)$:
$(D^{(n)}f)(a)=\sum_{b=0}^{\infty}D^{(n)}(a,b) f(b)$. In the natural
$\ell^{2}(\mathbb{Z}_{+})$ basis $\{\delta_{i}\}_{i=0}^{\infty}$ the
matrix elements are
\begin{equation}
  \label{eq:54}
  D^{(n)}(i,j)=\left\{
    \begin{array}{lll}
      B(i)\sqrt{\frac{W(i)}{W(i+1)}},& j=i+1, & i,j\leq L,\\
      C(i)\sqrt{\frac{W(i)}{W(i-1)}},& j=i-1, & i,j\leq L,\\
      1, & i=j, & i,j>L\\
      0 &\; &\text{otherwise}.
    \end{array}\right.
\end{equation}
Clearly, the operator $D^{(n)}$ is self-adjoint. Therefore by using
Theorem VIII.25 in \cite{reed1981}, similarly to Hahn ensemble in
\cite{gorin2008nonintersecting}, we get the convergence to the limit
operator $D$, defined by the corresponding three-diagonal Jacobi
matrix, in the strong resolvent sense.

It is easy to check that in the limit $N,a\to\infty, a\sim N$ the
ratio of the coefficients on the left hand side of equation
\eqref{eq:29} converges to 1
\begin{multline}
  \label{eq:31}
  \dfrac{B(a)}{C(a)}\sqrt{\dfrac{W(a-1)}{W(a+1)}}=\\
  =-\dfrac{1-q^{a-N}}{1-q^a}p^{-1}
  \sqrt{\dfrac{(1-q^a)(1-q^{a+1})}{(1-q^{N-a})(1-q^{N+1-a})}p^2
  	q^{2(N-a)-1}}
  \xrightarrow[N,a \to \infty]{} 1.
\end{multline}
Therefore for the ease of computation we can rewrite the difference
equation in the form
\begin{equation}
  \label{eq:32}
  \frac{B(a)}{C(a)}\sqrt{\frac{W(a-1)}{W(a+1)}}\kappa_m(a+1)+\kappa_m(a-1)=
  \sqrt{\dfrac{W(a+1)}{W(a)}}
  \left(\dfrac{A(m)}{B(a)}+1+\dfrac{C(a)}{B(a)}\right)\kappa_m(a).
\end{equation}

Then eigenvalues on the right hand side are:
\begin{equation}
  \label{eq:33}
  q^{\frac{1}{2}(a-N)}p^{-\frac{1}{2}}
  \sqrt{\dfrac{\left[N-a\right]_{q}}{\left[a+1\right]_{q}}}
  \left(\frac{q^{-m}(1-q^{m})(1+pq^m)}{1-q^{a-N}}+1-\frac{p(1-q^a)}{(1-q^{a-N})}\right),
\end{equation}
where $m=0,\dots,N$.
Substituting $p=q^{1-2n}$, $N=n+k-1$, $q=e^{-\gamma\frac{1}{n}}$, $a=nx$ and
taking the limit $n,k\to\infty$, we see that eigenvalues fill the
interval 
\begin{multline}
  \label{eq:34}
  \left[e^{\frac{\gamma}{2}(c+1-x)-\gamma}
    \sqrt{\frac{1-e^{-\gamma(c+1-x)}}{1-e^{-\gamma x}}}
    \left(\frac{(e^{(c+1)\gamma}+1)(e^{(c+1)\gamma}-1)}{1-e^{\gamma(c+1-x)}}+1-
      \frac{e^{2\gamma}(1-e^{-\gamma x})}{1-e^{\gamma(c+1-x)}}\right),
  \right.\\\left.
    e^{\frac{\gamma}{2}(c+1-x)-\gamma}
    \sqrt{\frac{1-e^{-\gamma(c+1-x)}}{1-e^{-\gamma x}}}
    \left(1-
      \frac{e^{2\gamma}(1-e^{-\gamma x})}{1-e^{\gamma(c+1-x)}}\right)
  \right].
\end{multline}
The corresponding limit operator $\widetilde{D}$ then acts as a
difference operator 
\begin{equation}
	\widetilde{D}f(x)=f(x+1)+f(x-1).
\end{equation}
The operator
$\widetilde{D}$ is self-adjoint and has simple purely continuous Lebesgues spectrum.
The correlation kernel $\mathbf{K}_{n}(a,b)$ is the projection to the
part of the spectrum that in the limit becomes
\begin{multline}
  \label{eq:35}
    \left[e^{\frac{\gamma}{2}(c+1-x)-\gamma}
    \sqrt{\frac{1-e^{-\gamma(c+1-x)}}{1-e^{-\gamma x}}}
    \left(\frac{(e^\gamma+1)(e^\gamma-1)}{1-e^{\gamma(c+1-x)}}+1-
      \frac{e^{2\gamma}(1-e^{-\gamma x})}{1-e^{\gamma(c+1-x)}}\right),
  \right.\\\left.
    e^{\frac{\gamma}{2}(c+1-x)-\gamma}
    \sqrt{\frac{1-e^{-\gamma(c+1-x)}}{1-e^{-\gamma x}}}
    \left(1-
      \frac{e^{2\gamma}(1-e^{-\gamma x})}{1-e^{\gamma(c+1-x)}}\right)
  \right].
\end{multline}
This spectral projection is given by the discrete sine kernel
\begin{equation}
	\mathbf{K}^{sine}_{\varphi}(u,v)=\frac{\sin(\varphi(u-v))}{\pi(u-v)},
\end{equation}
 since
the Fourier transform of the difference operator $\widetilde{D}$ from
$\ell^{2}(\mathbb{Z})$ to $L^{2}$ on unit circle $|z|=1$  is the
multiplication by the function $z+\bar{z}=2\Re z$ and has purely
continuous (double) spectrum $[-2,2]$. The maximum of the interval
\eqref{eq:35} is always $\geq 2$ while the minimum is inside of
$[-2,2]$. So in the $L^{2}$ space the spectral projection is the
multiplication by the characteristic function of the arc from
$e^{-i\varphi}$ to $e^{i\varphi}$. In the $\ell^{2}$ realization, this
is the integral operator with the discrete sine kernel with $\varphi$
given by the formula
\begin{equation}
  \label{eq:36}
  \varphi=\arccos\left(\frac{1}{2}e^{\frac{\gamma}{2}(c+1-x)-\gamma}
    \sqrt{\frac{1-e^{-\gamma(c+1-x)}}{1-e^{-\gamma x}}}
    \left(\frac{(e^\gamma+1)(e^\gamma-1)}{1-e^{\gamma(c+1-x)}}+
      1-\frac{e^{2\gamma}(1-e^{-\gamma x})}{1-e^{\gamma(c+1-x)}}\right)\right).
\end{equation}
Now as we are interested in the one-point correlation function, we
take limit $u\to v$ in the correlation kernel and get
$\rho(x)=\frac{\varphi}{\pi}$, which after a simplification becomes
\eqref{eq:limit-shape}. Similarly, taking $N=n+k-1$ and $p=q^{2-2n-k}$
and weight \eqref{eq:22}, we recover the spectral interval and the
limit density \eqref{eq:limit-shape-2} for the measure~\eqref{eq:5}. 

Since the points $a_{i}$ of the discrete \qKrawtchouk ensemble
correspond to the intervals where the upper boundary $f_{n}$ decays,
the density of these points is connected to the derivative of $f_{n}$
by the formula $f_{n}'(x)=1-2\rho^{(n)}_{1}(x)$. Thus from the limit density
we recover the limit shape by the formula \eqref{eq:20}. 
%%  To complete the proof we need to establish that the difference
%%  operators $D^{(n)}$, the operator $D$ and its closure $\overline{D}$
%%  are self-adjoint, bounded and have simple spectrum on the
%%  corresponding Hilbert space. All these questions are reduced to the
%%  Hambuerger moment problem along the lines of \cite[Lemma 3.1 and
%%  Theorem 4.1]{borodin2007asymptotics}. Hamburger moment problem for
%%  \qKrawtchouk ensemble is determinate for finite $n$ (see PhD thesis
%%  \cite{christiansen2004indeterminate} for further details), the limit
%%  $n\to\infty$ can be computed, but we omit this computation as not
%%  particularly interesting. 

\section{Fluctuations}
\label{sec:fluctuations}
In this section we prove Theorem \ref{central-limit-theorem}
applying the general approach of Breuer and Duits~\cite{breuer2017central}. Due to their result, it is
sufficient to establish the convergence of the coefficients in the
three-term recurrence relation for the corresponding orthogonal polynomials. Then we can apply the following theorem.
\begin{Th}[Theorem 2.5 in \cite{breuer2017central}]
  \label{thm:breuer-clt}
  Let $\{p^{(n)}_{m}(x)\}_{m=0}^{n-1}$ be normalized orthogonal
  polynomials of the polynomial ensemble $\mathcal{P}_{n}$ 
  satisfying the three-term recurrence relation
  \begin{equation*}
    \label{eq:38}
    xp^{(n)}_{m}(x)=a^{(n)}_{m+1}p^{(n)}_{m+1}(x)+b^{(n)}_{m+1}p^{(n)}_{m}(x)+a^{(n)}_{m}p^{(n)}_{m-1}(x),
  \end{equation*}
  and assume that there exists a subsequence $\{n_{j}\}_{j}$ and
  $a>0, b\in\mathbb{R}$ such that for any $k\in\mathbb{Z}$ we have
  \begin{equation*}
    \label{eq:39}
    a^{(n_{j})}_{n_{j}+k}\to a, \quad b^{(n_{j})}_{n_{j}+k}\to b, 
  \end{equation*}
  as $j\to\infty$. Then for any real-valued $f\in C^{1}(\mathbb{R})$
  we have
  \begin{equation*}
    \label{eq:40}
    X_{f}^{\mathcal{P}_{n_{j}}}-\mathbb{E}X_{f}^{\mathcal{P}_{n_{j}}} \to \mathcal{N}\left(0,\sum_{l\geq 1}l|\widehat{f}_{l}|^{2}\right), \quad \text{as}\; j\to\infty,
  \end{equation*}
  in distribution, where the coefficients $\widehat{f}_{l}$ are defined as
  \begin{equation*}
    \label{eq:41}
    \widehat{f}_{l}=\frac{1}{2\pi}\int_{0}^{2\pi} f(2a\cos \theta +b) e^{-il\theta}d\theta,
  \end{equation*}
  for $l\geq 1$. When $n_{j}=j$, that is the subsequence is the whole
  sequence, \eqref{eq:39} is equivalent to
  \begin{equation*}
    \label{eq:42}
    a^{(n)}_{n}\to a,\quad b^{(n)}_{n}\to b.
  \end{equation*}
\end{Th}
Similarly to the study of the lozenge tilings and Hahn polynomial
ensemble in \cite[Section 6.2]{breuer2017central}, we establish the required convergence for the normalized \qKrawtchouk polynomials. We start with the recurrence relation for the monic \qKrawtchouk polynomials $P_{n}(q^{-x})$ \cite[formula 14.15.4]{koekoek2010hypergeometric}:
\begin{equation*}
  \label{eq:43}
  q^{-x}
  P_{n}(q^{-x})=P_{n+1}(q^{-x})+
  \left[1-(A_{n}+C_{n})\right]P_{n}(q^{-x})+A_{n-1}C_{n}P_{n-1}(q^{-x}), 
\end{equation*}
where $P_{n}(q^{-x})=\frac{(q^{-N};q)_{n}}{(-pq^{n};q)_{n}}K_{n}(q^{-x};p,N;q)$, $q$-Pochhammer symbols are defined by the formula \eqref{eq:45} and
\begin{eqnarray} \label{eq:44-1}
  A_{n}=\frac{(1-q^{n-N})(1+pq^{n})}{(1+pq^{2n})(1+pq^{2n+1})},\\
  C_{n}=-pq^{2n-N-1}\frac{(1+pq^{n+N})(1-q^{n})}{(1+pq^{2n-1})(1+pq^{2n})}.
\end{eqnarray}
If the monic polynomials $P_{n}(x)$ satisfy the recursion relation
$$x P_{n}(x)=P_{n+1}(x)+\alpha_{n}P_{n}(x)+\beta_{n}P_{n-1},$$
then the normalized polynomials $p_{n}(x)$ satisfy the recursion
relation \eqref{eq:38} with $a_{n}=\sqrt{\beta_{n}}$,
$b_{n+1}=\alpha_{n}$ \cite{ismail2020encyclopedia}. Therefore for the
normalized \qKrawtchouk polynomials we have
\begin{eqnarray}
  a_{n}=\sqrt{A_{n-1}C_{n}},\\
  b_{n+1}=1-A_{n}-C_{n}. \label{eq:46-2}
\end{eqnarray}
Substituting  $p=q^{1-2n}, N=n+k-1, q=e^{-\gamma\frac{1}{n}}$ and
taking the limit $n\to\infty$, we get
\begin{equation}
  \label{eq:47}
  \begin{array}{l}
  a=\lim_{n\to\infty}a_{n}=
  \frac{1}{4}\sqrt{(e^{\gamma}-1)(e^{\gamma}+1)(1-e^{c\gamma})(1+e^{c\gamma})},\\
  b=\lim_{n\to\infty}b_{n}=\frac{1+e^{\gamma(c+1)}}{2}.    
  \end{array}
\end{equation}
The interval $[b-2a,b+2a]$ is exactly the support of the limit density
\eqref{eq:limit-shape}, written in terms of the variable $e^{t}$.
This can be easily verified by solving the equation
\begin{equation*}
  \label{eq:48}
  \mathrm{sgn}(-\gamma)\dfrac{e^{\gamma-\frac{\gamma  t}{2}}}{2}
  \dfrac{1-e^{\gamma(c-1)}}{\sqrt{(1-e^{\gamma
        t})(1-e^{\gamma(c+1-t)})}}=\pm 1.
\end{equation*}

Now we apply Theorem \ref{thm:breuer-clt} to get the desired result
for the probability (\ref{eq:4}). 

Similarly, substituting  $N=n+k-1$, $q=e^{-\gamma\frac{1}{n}}$ and
$p=q^{2-2n-k}$ into \eqref{eq:44-1}--\eqref{eq:46-2}, we obtain
\begin{equation}
  \label{eq:49}
  \begin{array}{l}
  a=\lim_{n\to\infty}a_{n}=
  \frac{e^{c\gamma}}{(1+e^{c\gamma})^{2}}
  \sqrt{2(e^{\gamma}-1)(e^{c\gamma}-1)(1-e^{c\gamma})(1+e^{\gamma(c+1)})},\\
  b=\lim_{n\to\infty}b_{n}=
  \frac{3e^{\gamma(c+2)}-e^{\gamma(c+1)}+3 e^{c\gamma}-e^{2c\gamma}}{(1+e^{c\gamma})^{2}}.      
  \end{array}
\end{equation}
The interval $[b-2a,b+2a]$ is exactly the support of the limit density
\eqref{eq:limit-shape-2}, written in terms of the variable $e^{t}$.
This can be easily verified by solving the equation
\begin{equation*}
  \label{eq:50}
  \mathrm{sign}(-\gamma)
  \frac{e^{\frac{\gamma}{2}(t-c)}}{2}
  \dfrac{1-e^{\gamma c}-e^{\gamma (c-t)}+e^{\gamma (c+1-t)}}
  {\sqrt{(1-e^{\gamma t})(1-e^{\gamma(c+1-t)})}}=\pm 1.  
\end{equation*}
Again, we apply Theorem \ref{thm:breuer-clt}  to finish the proof of
Theorem \ref{central-limit-theorem}.     

\section*{Conclusion and outlook}
\label{sec:conclusion-outlook}

\begin{figure}[!bth]
  \centering
  \includegraphics[width=12cm]
  {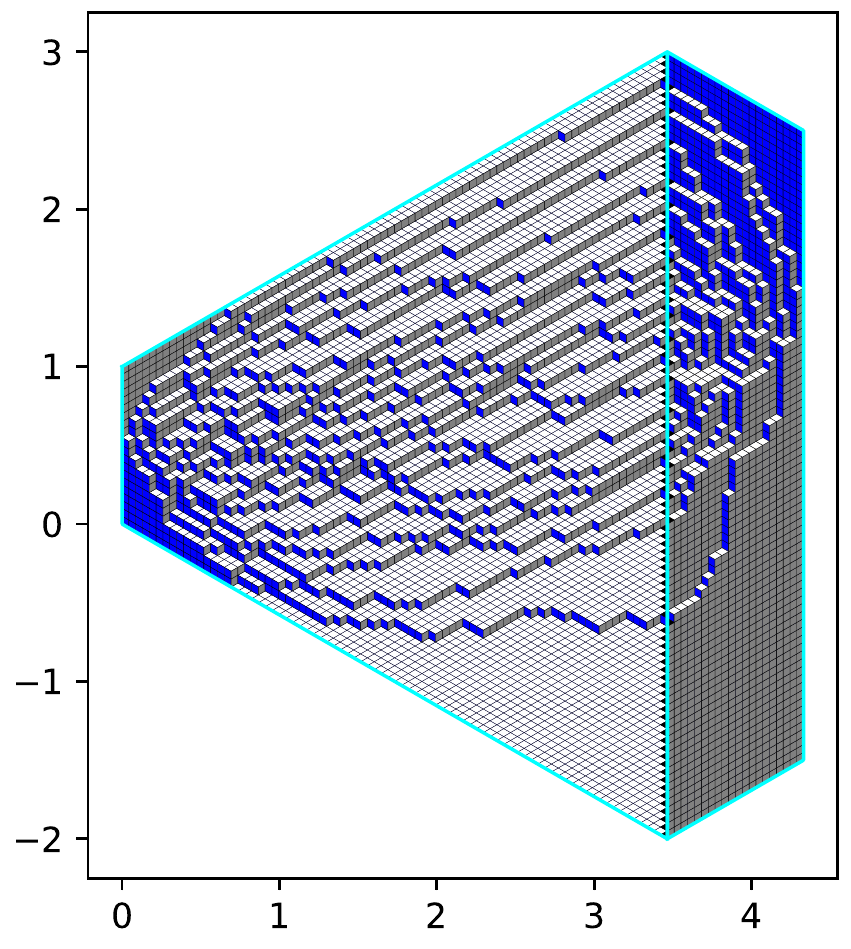}
  \caption{Lozenge tiling of skew glued hexagon, Young diagrams
    $\lambda,\bar\lambda'$ correspond to the positions of right (left)
  triangles on the main diagonal.}
  \label{fig:lozenge-tiling}  

\end{figure}

In the present paper we have discussed the \qKrawtchouk polynomial
ensemble from the point of view of the skew Howe duality. We have
connected the limit densities of the point processes to the limit
shapes of random Young diagrams that parameterize the decomposition of
the exterior algebra
$\bigwedge\left(\mathbb{C}^{n}\otimes\left(\mathbb{C}^{k}\right)^*\right)$
into the irreducible $\GL_{n}\times\GL_{k}$-modules with respect to
the principal-principal and principal-inverse principal
specializations of the character measure \eqref{eq:3}. As was
demonstrated in the paper \cite{nazarov2021skew}, such measures can be
thought of as the measures induced on the main diagonal of the lozenge
tiling for the skew-glued hexagon from the tilings weighted by
$q^{\text{Volume of boxes}}$. Volume-weighted tilings of the hexagon
can be related to the symmetric $(\GL_{n},\GL_{k})$ Howe duality (See
\cite{cohn1998shape,okounkov2001infinite,okounkov2003correlation,gorin2021lectures}). It remains an interesting question to describe the limit surface of such volume-weighted skew-glued tilings that represent a skew Howe duality counterpart to the ordinary volume-weighted tilings of the hexagon.

In Fig.~\ref{fig:lozenge-tiling} we present a random lozenge tiling
of a skew glued hexagon that contributes to a random Young diagram
$\lambda$ chosen from $\mu_{n,k}(\lambda|q)$, with ${n=20}$, ${k=80}$, ${\gamma=-0.5}$, that corresponds to the positions of the right triangles on the main diagonal. The complement conjugate diagram $\bar\lambda'$ corresponds to the positions of the left triangles on the main diagonal.

\section*{Acknowledgements}
This work is supported by the Russian Science Foundation under grant
No.~21-11-00141. We thank Travis Scrimshaw for the valuable comments on the draft version of the text.

\printbibliography

\end{document}